\documentclass[12pt]{article}
\usepackage{amsfonts,amssymb}

\def\C{\mathbb{C}}

\def\bq{ \begin{equation} }
\def\eq{ \end{equation} }
\def\ben{ \begin{eqnarray} }
\def\en{ \end{eqnarray} }

\def\frac#1#2{{#1\over #2}}

\def\on#1#2{\mathop{\vbox{\ialign{##\crcr\noalign{\kern2pt}
$\scriptstyle{#2}$\crcr\noalign{\kern2pt\nointerlineskip}
\kern-2pt$\hfil\displaystyle{#1}\hfil$\crcr}}}\limits}

\font\Sets=msbm10

\def\Integ {\hbox{\Sets Z}}
\begin{document}

\baselineskip=15pt
\vspace{1cm} \noindent {\LARGE \textbf {Integrable matrix
 equations related to pairs of compatible associative algebras}}
\vskip1cm \hfill
\begin{minipage}{13.5cm}
\baselineskip=15pt {\bf
 A V Odesskii ${}^{1,} {}^{2}$ and
 V V Sokolov ${}^{1}$} \\ [2ex]
{\footnotesize ${}^1$ Landau Institute for Theoretical Physics,
Moscow, Russia
\\
${}^{2}$  School of Mathematics, The University of
Manchester, UK}\\
\vskip1cm{\bf Abstract}

We study associative multiplications
in semi-simple associative algebras over $\C$ compatible with the
usual one. An interesting class of such multiplications is related to
the affine Dynkin diagrams of $\tilde A_{2 k-1},$ $\tilde D_{k},$
$\tilde E_{6},$ $\tilde E_{7},$ and $\tilde E_{8}$-type.
In this paper we investigate in details the
multiplications of the $\tilde A_{2 k-1}$-type and integrable
matrix ODEs and PDEs generated by them.

\end{minipage}

\vglue1cm \textbf{Address}: Landau Institute for Theoretical
Physics, Kosygina 2, 119334, Moscow, Russia

\textbf{E-mail}: odesskii@itp.ac.ru, \, sokolov@itp.ac.ru \newpage

\centerline{\Large\bf Introduction}
\medskip

In the papers \cite{golsok1,golsok2,golsok3,golsok4} different
applications of the notion of compatible Lie brackets to the
integrability theory have been considered.  Various odd examples
of compatible Lie structures have been presented in \cite{bolsbor,
golsok1, golsok2}.

A pair of compatible associative multiplications is more rigid
algebraic structure then a pair of compatible Lie brackets and
therefore there is a chance to construct more developed theory and
some series of important examples of such multiplications.

In \cite{odsok} we have studied multiplications  compatible with the
matrix multiplication or, in other words, linear deformations of the
matrix product. It turns out that these deformations of the matrix
algebra are in one-to-one correspondence with representations of
certain algebraic structures, which we call $M$-structures. The case
of direct sum of several matrix algebras corresponds to
representations of the so-called $PM$-structures (see \cite{odsok}).

The main result of \cite{odsok} is a description of an important
class of $M$ and $PM$-structures. The classification of these
structures naturally leads to the Cartan matrices of affine Dynkin
diagrams of the $\tilde A_{2 k-1},$ $\tilde D_{k},$ $\tilde
E_{6},$ $\tilde E_{7},$ and $\tilde E_{8}$-type.

It this paper, we investigate  integrable equations related to the
affine Dynkin diagrams of $\tilde A$-type in details. Recall some
results of the papers \cite{golsok4,odsok,golsok1} adopted for our
goals. For simplicity we formulate these results for the matrix
algebra $Mat_n$ but all of them remain to be true for any
semi-simple associative algebra.

A multiplication $\circ$ defined on the vector space $Mat_n$ of
all $n\times n$ matrices is said to be compatible with the matrix
multiplication if the product
\begin{equation} \label{pensil}
X \bullet Y= X  Y+ \lambda \,X\circ Y
\end{equation}
is associative for any constant $\lambda$.

Since the matrix algebra is rigid, the multiplication
(\ref{pensil}) is isomorphic to the matrix multiplication for
almost all values of the parameter $\lambda$. This means that
there exists a formal series of the form
\begin{equation}
S_{\lambda}={\bf 1}+ R \ \lambda+ O(\lambda^2),    \label{RS}
\end{equation}
with the coefficients being linear operators on $Mat_n$, such that
\begin{equation}
S_{\lambda}(X) \, S_{\lambda}(Y)= S_{\lambda} \Big( X Y+\lambda \
X \circ Y\Big). \label{sog2}
\end{equation} It follows from this formula that the multiplication $\circ
$ is given by
\begin{equation} \label{mult2}
X \circ Y =R(X) Y+X R(Y)-R(X Y).
\end{equation}
 The series $S_\lambda$
is defined up to a transformation $S_{\lambda}\rightarrow
T_{\lambda}\, S_{\lambda}\,T_{\lambda}^{-1},$ where $T_{\lambda}$ is
an arbitrary element of $Mat_n$ of the form $T_{\lambda}={\bf 1}+t
\lambda+O(\lambda^2)$. This leads to the transformation
\begin{equation} \label{ad}
R \rightarrow R + ad_t,
\end{equation}
which does not change the multiplication (\ref{mult2}).

In the matrix case it is convenient to write down the operator $R$
in the form
\begin{equation}\label{Rmat}
R(x)=a_1 \,x \,b^1+...+a_p \,x\, b^p+c\, x, \qquad \qquad a_i,
b^i, c \in Mat_n
\end{equation}
with $p$ being smallest possible.

A general construction from \cite{golsok4} establishes a
relationship between pairs of compatible associative
multiplications and integrable top-like systems. Let $\circ$ be a
multiplication (\ref{mult2}) compatible with the matrix product.
Consider the following matrix differential equation
\begin{equation}\label{Emat}
\frac{d x}{dt}=[R(x)+R^*(x), x],
\end{equation}
where $*$ stands for the adjoint operator with respect to the
bi-linear form $\mbox{trace}\,(x y).$ This means that, if $Q$ is a
linear operator in $Mat_n$, then $Q^*$ is defined by the equation
$\mbox{trace}(Q(x)y)=\mbox{trace}(xQ^*(y))$ for all $x,y\in Mat_n$.
In particular,
$$
R^*(x)=b^1 \,x \,a_1+...+ b^p\,x\,a_p + x \, c.
$$

It turns out \cite{golsok4} that equation (\ref{Emat}) possesses
the Lax representation
$$
\frac{d L}{dt}=[A,\,L],
$$
where
\begin{equation}\label{AL}
L=\left(S_{\lambda}^{-1}\right)^* (x),  \qquad A=\frac{1}{\lambda}
S_{\lambda}(x).
\end{equation}
To make these formulas constructive, we should find $S_{\lambda}$
and $S_{\lambda}^{-1}$ in a closed form i.e. as analytic
operator-valued functions. As usual, the integrals of motion for
(\ref{Emat}) are given by coefficients of different powers of
$\lambda$ in $\mbox{trace}\,(L^k), \,\, k=1,2...$

The simplest example of a multiplication compatible with
the matrix product corresponds to the Dynkin diagram  $\tilde A_{1}.$
This multiplication $\circ$  is given by $x\circ y=x c y,$ where $c$ is an
arbitrary matrix. In this case we have
$$
R(x)=c x, \qquad S_{\lambda} (x)=(1+\lambda c)\, x, \qquad L=x
(1+\lambda c)^{-1}.
$$
The Lax equation is equivalent to the  well known integrable matrix
equation
\begin{equation}\label{xc}
\frac{d x}{dt}=x^2\,c-c\, x^2.
\end{equation}
The simplest integrals for (\ref{xc}) are given by $H_{k,0}
=\mbox{trace}\,(x^k)$. Moreover, the $L$-operator produces  an
infinite set of homogeneous integrals $H_{i,j}$, where $i$ and $j$
are degrees with respect to $x$ and $c$.  For example,
\begin{equation}\label{H2}
 H_{1,1} =\mbox{trace}\,(x
c)\, , \quad H_{2,1} =\mbox{trace}\,(x^2 c)\, , \quad H_{2,2} =
\mbox{trace}\,(2c^2 x^2+ c x c x).
\end{equation}
Equation (\ref{xc}) is Hamiltonian one with respect to the
standard matrix linear Poisson bracket given by the Hamiltonian
operator $ad_x$ and Hamiltonian function $H_{2,1}$. Indeed,
(\ref{xc}) can be rewritten in the form (\ref{Emat}) as
$x_{t}=[x,\,x c+c x]$.

Matrix equations of arbitrary size like (\ref{xc}) are important
because of possibility to make different reductions. For the most
trivial reduction one may regard $x$ as a block-matrix. In this
case (\ref{xc}) becomes a system of several matrix equations for
the block entries of $x$. Under reduction $x^T=-x, \, c^T=c$
(\ref{xc}) is a commuting flow for the $n$-dimensional Euler
equation \cite{man, fed, miksok}. Another reduction of (\ref{xc})
was mentioned in \cite{miksok}. Let $x$ and $c$ in (\ref{xc}) be
represented by matrices of the form
\begin{equation}\label{MCmatr}
x=\left( \begin{array}{cccccc}
0&u_1&0&0&\cdot&0\\
0&0&u_2&0&\cdot&0\\
\cdot&\cdot&\cdot&\cdot&\cdot&\cdot\\
0&0&0&0&\cdot&u_{N-1}\\
u_{N}&0&0&0&\cdot&0
\end{array}\right)\, , \quad
c=\left( \begin{array}{cccccc}
0&0&0&\cdot&0&J_N\\
J_1&0&0&\cdot&0&0\\
0&J_2&0&\cdot&0&0\\
\cdot&\cdot&\cdot&\cdot&\cdot&\cdot\\
0&0&0&\cdot&J_{N-1}&0
\end{array}\right)\, ,
\end{equation}
where $u_k$ and $J_k$ are block matrices (of any dimension).  It
follows from equation (\ref{xc}) that $u_k$ satisfy the
non-abelian Volterra equation
$$
\frac{d}{dt}u_k=u_k \circ u_{k+1}\circ J_{k+1}-J_{k-1}\circ
u_{k-1}\circ u_k, \quad k\in\Integ _{N}\, .
$$

Integrable equation (\ref{xc}) contains one arbitrary constant matrix $c$. Such equations
have been systematically investigated in \cite{miksok}. In
Sections 1,2 we present a series of integrable matrix equations,
which depend on two constant matrices related by certain algebraic
relations providing the integrability. In Section 3 these results
are generalized to the case of several matrix unknowns.
Skew-symmetric reductions as well as reductions of
(\ref{MCmatr})-type are available for our matrix ODEs.

Different applications of compatible associative multiplications are
related to integrable deformations of the principle $GL(n)$-chiral
model (see \cite{golsok1}). Let $\circ$ be associative
multiplication compatible with the  $Mat_n\oplus Mat_n$ product,
$S_{\lambda}$ be the series (\ref{RS}) with the property
(\ref{sog2}). Define operators $T_1, T_2$ on $Mat_n$ by means of the
following decompositions
\begin{equation} \label{decom}
S_{\lambda}(u,0)=(u,0)+\lambda (\cdot, \, T_1(u))+O(\lambda^2),
\qquad S_{\lambda}(0,v)=(0,v)+\lambda (T_2(v),\,
\cdot)+O(\lambda^2).
\end{equation}
Then the following hyperbolic system
\begin{equation} \label{chir}
u_t=[u,\, T_2(v)], \qquad v_{\tau}=[v,\, T_1(u)], \qquad \quad
u,v\in Mat_n
\end{equation}
possesses a zero-curvature representation $[L,M]=0$, where
$$
L=\frac{d}{dy}+\frac{1}{\lambda} S_{\lambda}(u,0), \qquad
M=\frac{d}{dx}+\frac{1}{\lambda} S_{\lambda}(0,v).
$$
If $T_1=T_2={\bf 1}$, then (\ref{chir}) is just the principle chiral
model \cite{pohl}. In Section 3 we describe a class of integrable
models (\ref{chir}) related to $PM$-structures of the $\tilde A_{2
k-1}$-type.

\section{Associative product of the $\tilde A_{3}$-type.}

The $M$-structure related to $\tilde A_{3}$ (see \cite{odsok} ) is
defined by two arbitrary constant matrices $A$ and $B$  such that
\begin{equation}\label{squ}A^2=B^2={\bf 1}.\end{equation}
The corresponding associative multiplication is given by
(\ref{mult2}), where
\begin{equation}\label{Rab}
R(x)=A x B+ B A x.
\end{equation}
This structure leads to the following integrable matrix equation
\begin{equation}\label{abx}
x_{t}=[x,\,\, B x A+A x B+x B A+B A x].
\end{equation}
The Lax representation (\ref{AL}) for (\ref{abx}) is given by the
following  explicit formulas for $S_{\lambda}$ and
$S_{\lambda}^{-1}$:
$$
S_{\lambda}(x)=\frac{1-q}{2} \,B x B+\frac{1+q}{2}\, x+\lambda (A
x B+B A x),
$$
$$
S_{\lambda}^{-1}(x)=\frac{1}{q} ({\bf 1}+\lambda
K)^{-1}\Big(\frac{q-1}{2} \,B x B+\frac{1+q}{2}\, x+\lambda (A B
x-A x B)\Big),
$$
where $q=\sqrt{1-4\lambda^2},\,\, $ $K=A B+B A$. Note that both $A$
and $B$ commute with $K$.

The canonical form (\ref{Rab}) for the operator $R$ with respect
to transformations (\ref{ad}) is homogeneous in $A$ and $B$ that
gives rise to short expressions for $S_{\lambda}$ and
$S_{\lambda}^{-1}$ irrationally depending on $\lambda$. A
different form
\begin{equation}\label{Rab1}
\bar R(x)=A x B+ B A x+B x-x B
\end{equation}
provides the simplest form of the operator (\ref{RS}) with the
property (\ref{sog2}):
\begin{equation}\label{SiS}
\bar S_\lambda(x)=x+ \lambda \bar R(x)
\end{equation}
but the expression for $\bar S_{\lambda}^{-1}$ is more
complicated:
$$
\bar S_{\lambda}^{-1}(x)=\frac{1}{4 \lambda^2-1} ({\bf 1}+\lambda
K)^{-1}\Big(\lambda \,(A x B+B x+B A x)+2 \lambda^2 \,(B x B+B A x
B+A x)$$$$-(\lambda + 2 \lambda^2 K)\, x B-(1+\lambda K-2
\lambda^2) \, x\Big).
$$
The corresponding Lax operators (\ref{AL}) are rational in
$\lambda$.

The simplest linear and quadratic first integrals for (\ref{abx})
generated by the $L$-operator are given by
$$
 H_{1,1} =\mbox{trace}\,[x
(A B+B A)], \quad H_{1,2} =\mbox{trace}\,[x (A B A B+B A B A)],
\quad H_{2,1} =\mbox{trace}\,[B A x^2+A x B x]\, ,
$$
and
$$
H_{2,2} =\mbox{trace}\,[2 B A B A x^2+2 A B A x B x+ 2 B A B x A x+
A B x A B x+B A x B A x].
$$

{\bf Remark.} Equation (\ref{abx}) is the simplest example of matrix
integrable equation with constant matrices related by constraints.
If we assume that $A$, $B$, $x$ are block-matrices of the form:

$$A=\pmatrix{{\bf 1}&0&\cr 0&-{\bf 1}\cr}, \qquad B=\pmatrix{P& {\bf 1}+P&\cr {\bf
1}-P&-P\cr},\qquad x=\pmatrix{x_{1,1}&x_{1,2}&\cr
x_{2,1}&x_{2,2}\cr},$$ then we obtain an integrable  ODE with one
arbitrary constant matrix $P$ and four matrix unknowns $x_{i,j}$.
Furthermore, assuming that $P$ is scalar matrix, we obtain an
integrable matrix ODE with constant coefficients.

Equation (\ref{abx}) admits the following skew-symmetric reduction
\begin{equation}
\label{red} x^T=-x, \qquad B=A^T.
\end{equation}
Different integrable $so(n)$-models provided by reduction
(\ref{red}) are in one-to-one correspondence with equivalence
classes with respect to the $SO(n)$ gauge action of $n\times n$
matrices $A$ such that $A^2={\bf 1}$. For the real matrix $A$, a
canonical form for such equivalence class can be chosen as
\begin{equation}\label{Mat}
A=\pmatrix{{\bf 1}_p&T&\cr 0&{\bf 1}_{n-p}\cr}.
\end{equation}
Here ${\bf 1}_s$ stands for the unity $s\times s$ matrix and
$T=(t_{ij}),$ where $t_{ij}=\delta_{ij} \alpha_i$. This canonical
form is defined by the discrete natural parameter $p$ and continuous
parameters $\alpha_1,\dots, \alpha_r$, where $p\le n/2,\,\,$
$r=\min(p,n-p)$. A different canonical form for $A$ will be given in
the next section.

For example, in the case $n=4$ the equivalence classes with $p=2$
and $p=1$ give rise to the Steklov  and the Poincare integrable
models \cite{stek1,puan}, correspondingly. Thus, whereas
(\ref{xc}) is a matrix version of the $so(4)$ Schottky-Manakov top
\cite{shot,man}, equation (\ref{abx})-(\ref{Mat}) with $p=[n/2]$
and $p=1$ can be regarded as $so(n)$ generalizations for the
$so(4)$ Steklov and Poincare models (for different generalization
see \cite{rat, bf}).

\section{Associative products of the $\tilde A_{2 k-1}$-type. Case of one
matrix variable.}

Let $A, B, C\in Mat_n$ be matrices that satisfy the following set of
relations
\begin{equation} \label{rel1}
A^{k}=B^{k}={\bf 1}, \end{equation}
\begin{equation}
\label{rel2}B^iA^j=\frac{\varepsilon^{-j}-1}{\varepsilon^{-i-j}-1}A^{i+j}+
\frac{\varepsilon^i-1}{\varepsilon^{i+j}-1}B^{i+j}, \qquad
i+j\ne0\quad \mbox{mod}\,\, k,
\end{equation}
\begin{equation} \label{rel3}B^iA^{k-i}={\bf 1}+(\varepsilon^i-1)C,
\end{equation}
where $\varepsilon=\exp (2\pi i/k).$ It follows from results of
\cite{odsok} that formula (\ref{mult2}) with
$$
R(x)=\sum_{i=1}^{k-1} \frac{1}{\varepsilon^i-1} A^{k-i}\,x\, B^i+
C x,
$$
defines the associative multiplication $\circ $ compatible with the
standard matrix product in $Mat_n.$ According to \cite{odsok} this
multiplication corresponds to the affine Dynkin diagram of the
$\tilde A_{2 k-1}$-type.

Let $T$ be any matrix such that $T^k-{\bf 1}$ is invertible and
\begin{equation} \label{rel4}
A\, T=\varepsilon T \, A.
\end{equation}
Then
\begin{equation} \label{rel5}
B= (\varepsilon T-{\bf 1}) (T-{\bf 1})^{-1}\, A, \qquad C=T(T-{\bf
1})^{-1}
\end{equation}
satisfy identities (\ref{rel1})-(\ref{rel3}). This is a convenient
way to resolve this system of identities.

{\bf Proposition 1.} The operator $S_{\lambda}$ with the property
(\ref{sog2}) is defined by
\begin{equation} \label{ss1}
S_{\lambda} (x)= \sum_{i=0}^{k-1} P_i x B^i,\qquad
P_i=\frac{\lambda}{\varepsilon^i \mu -1} (\mu T-{\bf 1}) (T-{\bf
1})^{-1} A^{k-i},
\end{equation}
where
$$
\mu=\left(\frac{2+(k+1) \lambda}{2-(k-1) \lambda}\right)^{1/k}.
$$
The inverse operator is given by the formula
\begin{equation} \label{ss2}
S_{\lambda}^{-1} (x)= \sum_{i=0}^{k-1} Q_i x B^i,\qquad
Q_i=\frac{(\mu^k-1)^2}{\lambda k^2 \mu^{k-1}(\mu- \varepsilon^i)}
(\varepsilon^{k-i} T-{\bf 1}) (\varepsilon^{k-i} \mu T-{\bf
1})^{-1} A^{k-i}.
\end{equation}
The corresponding Lax operator (\ref{AL}) is given by
\begin{equation} \label{ss1} L(x)=\sum_{i=0}^{k-1} B^i x
Q_i.
\end{equation}

{\bf Remark.} If we take operator $R$ in the following equivalent
but different form
$$\bar R(x)=\sum_{i=1}^{k-1} \frac{1}{\varepsilon^i-1}
(A^{k-i}-1)\,x\, B^i+ (C+\sum_{i=1}^{k-1}
\frac{1}{\varepsilon^i-1}B^i) \,x,
$$
then the operator $S_{\lambda}$ has the simplest form (\ref{SiS}),
but the formula for $S_{\lambda}^{-1} (x)$ is more complicated.

Using Proposition 1, one can calculate the integrals of motion for
equation (\ref{Emat}) as coefficients of $\mbox{trace}\,(L^k),
\,\, k=1,2...$.  For example, the simplest quadratic integral is
given by
$$
H=\mbox{trace}\,(x R(x)).
$$
Equation (\ref{Emat}) can be written in the Hamiltonian form as
$$
\frac{d x}{dt}=[x,\, \mbox{grad}\, H].
$$

The skew-symmetric reduction is given by
$$B=(A^t)^{-1}, \qquad C=\frac{\varepsilon}{1-\varepsilon}(A^t
A-1),$$ where $A$ should satisfy the following additional
constraint
$$(A^t A-\varepsilon^{-1}) \, A \, (A^t A-1)=\varepsilon (A^t
A-1)\,A \, (A^t A-\varepsilon^{-1}).$$

It turns out that \begin{equation} \label{cred} A=\sqrt{z_1} \,
e_{k,1}+ \sum_{i=2}^k \sqrt{z_i} \, e_{i-1,i},\end{equation} where
$z_{j+1}=f(z_j),$
$$
f(z)=\frac{1}{1+\varepsilon-\varepsilon z},
$$
defines a skew-symmetric reduction with one arbitrary parameter
$z_1$. In particular, for $k=2$ we have
$$
A=\pmatrix{0&\alpha&\cr \alpha^{-1}&0\cr},
$$
where $\alpha$ is an arbitrary parameter.

{\bf Remark.} Any block-diagonal matrix with blocks (of  different
side) of the form (\ref{cred}) also provides a non-trivial
skew-symmetric reduction. If we take $p$ blocks of the size 1 and
the remaining blocks of the size 2, we get a canonical form for
matrix $A$ equivalent but different from (\ref{Mat}).

\section{Associative products of the $\tilde A_{2 k-1}$-type. Case of $m$
matrix variables.}

In this section we consider associative multiplications compatible
with the product in the direct sum of $m$ copies of $Mat_n$. We
use the notation $x=(x_1,...,x_m)$, where $x_1,...,x_m\in Mat_n,$
for $x\in (Mat_n)^m.$ The standard  multiplication in $(Mat_n)^m$
is given by the formula $xy=(x_1 y_1,...,x_m y_m).$ Using results
obtained in \cite{odsok}, the following statement can be proved.

{\bf Proposition 2.} Let $T,B$ be matrices such that $$B^k={\bf
1}, \qquad BT=\varepsilon TB.$$ Fix generic numbers
$\lambda_1,...\lambda_m,t_1,...,t_m\in\C$ such that
$T^k-\lambda_{\alpha}^k$ are invertible for all $\alpha=1,...,m.$
Then formula (\ref{mult2}) with
$$R(x)=(R_1(x),...,R_m(x)),$$
$$R_{\alpha}(x)=\sum_{(i,\beta)\ne(0,\alpha)}\frac{t_{\beta}}
{\varepsilon^i\lambda_{\beta}-\lambda_{\alpha}}\cdot\frac{T-\lambda_{\alpha}}
{\varepsilon^{-i}T-\lambda_{\beta}}B^{-i}x_{\beta}B^i+\Big(\frac{t_{\alpha}}{T-\lambda_
{\alpha}}-\sum_{(i,\beta)
\ne(0,\alpha)}\frac{\varepsilon^it_{\beta}}{\varepsilon^i\lambda_{\beta}-\lambda_{\alpha}}\Big)x_{\alpha}$$
defines the following associative product compatible with the
standard product in $(Mat_n)^m:$
$$x\circ y=((x\circ y)_1,...,(x\circ y)_m),$$
where
$$(x\circ y)_{\alpha}=
\sum_{(i,\beta)\ne(0,\alpha)}\frac{t_{\beta}}
{\varepsilon^i\lambda_{\beta}-\lambda_{\alpha}}\Big(\frac{T-\lambda_{\alpha}}
{\varepsilon^{-i}T-\lambda_{\beta}}B^{-i}x_{\beta}B^iy_{\alpha}+$$$$
x_{\alpha}\frac{T-\lambda_{\alpha}}
{\varepsilon^{-i}T-\lambda_{\beta}}B^{-i}y_{\beta}B^i-\frac{T-\lambda_{\alpha}}
{\varepsilon^{-i}T-\lambda_{\beta}}B^{-i}x_{\beta}y_{\beta}B^i\Big)+t_{\alpha}x_{\alpha}\frac{1}{T-\lambda_
{\alpha}}y_{\alpha}-\sum_{(i,\beta)
\ne(0,\alpha)}\frac{\varepsilon^it_{\beta}}{\varepsilon^i\lambda_{\beta}-\lambda_{\alpha}}x_{\alpha}y_{\alpha}.$$

{\bf Remark.} The formulas from Section 2 after elimination of
matrices $A^i$, $C$ via (\ref{rel5}) coincides with the
corresponding formulas from Proposition 2 with $m=1.$

{\bf Proposition 3.} For the multiplication from Proposition 2 the
operator (\ref{RS}) with the property (\ref{sog2}) is given by
$$S_{\lambda}(x)=(S_1(x),...,S_m(x)),$$
where
$$S_{\alpha}(x)=\sum_{i,\beta}\frac{\lambda t_{\beta}(T-\mu_{\alpha})}
{(\varepsilon^i\lambda_{\beta}-\mu_{\alpha})(\varepsilon^{-i}T-\lambda_{\beta})}B^{-i}x_{\beta}B^i.$$
The Lax operator (\ref{AL}) is defined by the formulas:
$$L(x)=(L_1(x),...,L_m(x)),$$ where

$$L_{\alpha}(x)=-\sum_{i,\beta}\frac{\prod_s(\lambda_{\alpha}^k-\mu_s^k)(\lambda_s^k-\mu_{\beta}^k)}
{k^2\lambda_{\alpha}^{k-1}\mu_{\beta}^{k-1}(\lambda_{\alpha}-
\varepsilon^i\mu_{\beta})\prod_{s\ne\alpha}(\lambda_{\alpha}^k-\lambda_s^k)\prod_{s\ne\beta}
(\mu_{\beta}^k-\mu_s^k)}\cdot B^i x_{\beta}
\frac{T-\lambda_{\alpha}}{t_{\alpha}\lambda
(\varepsilon^{-i}T-\mu_{\beta})}B^{-i}.$$ Here
$\mu_{\alpha}=\mu_{\alpha}(\lambda)$ are algebraic functions in
$\lambda$ such that $\mu_{\alpha}(0)=\lambda_{\alpha}$ and
$$\sum_{\gamma}\frac{t_{\gamma}\lambda_{\gamma}^{k-1}}{\lambda_{\gamma}^k-\mu_{\alpha}^k}=\frac{1}{k\lambda},
\qquad \alpha=1,...,m.$$

{\bf Remark.} In some sense, the Lax operator constructed in
Proposition 3 is a multi-pole analog  of operator (\ref{ss1}).

To find the explicit form of integrable equations related to this
Lax operator, one needs several first coefficients of the Taylor
expansions of $S_{\lambda}$ and $L.$ We have
$$\mu_{\alpha}=\lambda_{\alpha}-t_{\alpha}\lambda-t_{\alpha}\Big(\frac{k-1}{2}\cdot
\frac{t_{\alpha}}{\lambda_{\alpha}}+k\sum_{\gamma\ne\alpha}
\frac{t_{\gamma}\lambda_{\gamma}^{k-1}}{\lambda_{\gamma}^k-\lambda_{\alpha}^k}\Big)\lambda^2+O(\lambda^3),$$
$$S_{\alpha}(x)=x_{\alpha}+\Big(\frac{t_{\alpha}}{T-\lambda_{\alpha}}-\frac{k-1}{2}\cdot
\frac{t_{\alpha}}{\lambda_{\alpha}}-k\sum_{\gamma\ne\alpha}
\frac{t_{\gamma}\lambda_{\gamma}^{k-1}}{\lambda_{\gamma}^k-
\lambda_{\alpha}^k}\Big) x_{\alpha}\,\lambda+$$$$
\sum_{(i,\beta)\ne(0,\alpha)}\frac{t_{\beta}(T-\lambda_{\alpha})}
{(\varepsilon^i\lambda_{\beta}-\lambda_{\alpha})(\varepsilon^{-i}T-\lambda_{\beta})}
B^{-i}x_{\beta}B^i\,\lambda+O(\lambda^2),$$
$$L_{\alpha}(x)=x_{\alpha}-x_{\alpha}\Big(\frac{t_{\alpha}}{T-\lambda_{\alpha}}-\frac{k-1}{2}\cdot
\frac{t_{\alpha}}{\lambda_{\alpha}}-k\sum_{\gamma\ne\alpha}
\frac{t_{\gamma}\lambda_{\gamma}^{k-1}}{\lambda_{\gamma}^k-
\lambda_{\alpha}^k}\Big)\,\lambda$$
$$-\sum_{(i,\beta)\ne(0,\alpha)}B^ix_{\beta}\frac{t_{\beta}(T-\lambda_{\alpha})}
{(\varepsilon^i\lambda_{\beta}-\lambda_{\alpha})(\varepsilon^{-i}T-\lambda_{\beta})}B^{-i}\lambda+O(\lambda^2).$$

The integrable system (\ref{Emat}) has the form
$$\frac{d
x_{\alpha}}{dt}=\Big[\frac{t_{\alpha}}{T-\lambda_{\alpha}}x_{\alpha}+x_{\alpha}\frac{t_{\alpha}}{T-\lambda_{\alpha}}
+
$$
$$+\sum_{(i,\beta)\ne(0,\alpha)}\Big(\frac{t_{\beta}(T-\lambda_{\alpha})}
{(\varepsilon^i\lambda_{\beta}-\lambda_{\alpha})(\varepsilon^{-i}T-\lambda_{\beta})}B^{-i}x_{\beta}B^i+
B^ix_{\beta}\frac{t_{\beta}(T-\lambda_{\alpha})}
{(\varepsilon^i\lambda_{\beta}-\lambda_{\alpha})(\varepsilon^{-i}T-\lambda_{\beta})}B^{-i}\Big),\,\,x_{\alpha}\Big].$$
It can be written in the Hamiltonian form as
$$
\frac{d x}{dt}=[x,\, \mbox{grad}\, H],
$$
where $ H=\sum_{\alpha}\mbox{trace}\,(x_{\alpha} R_{\alpha}(x)).$

{\bf Remark.} All these constructions including the integrable
system are valid for the case $k=1$ also. In this case $B={\bf 1}$,
$\varepsilon=1$ and $T$ is an arbitrary matrix. If also $m=1$, then
the system is equivalent to (\ref{xc}).

Another application of Propositions 2,3 in the case  $m=2$ is the
integrable PDE (\ref{chir}). To write down this system explicitly,
we need to find the operators $T_1$ and $T_2$ defined by
(\ref{decom}). The asymptotic formula for $S_\lambda$ given above
yields
$$T_1(x)=\sum_i\frac{t_1(T-\lambda_2)}
{(\varepsilon^i\lambda_1-\lambda_2)(\varepsilon^{-i}T-\lambda_1)}B^{-i}xB^i,\qquad
T_2(x)=\sum_i\frac{t_2(T-\lambda_1)}
{(\varepsilon^i\lambda_2-\lambda_1)(\varepsilon^{-i}T-\lambda_2)}B^{-i}xB^i.$$

\section{Conclusion.}
In this paper we have described  representations of $M$ and
$PM$-structures related to the affine Dynkin diagrams of the
$\tilde A_{2 k-1}$-type and have presented Lax operators for the
corresponding matrix top-like systems. Actually, these systems are
bi-hamiltonian models \cite{magri}. We are planning to describe the Hamiltonian
properties of these systems in a separate paper. Besides of that
representations of $PM$-structures related to the affine Dynkin
diagrams of the $\tilde D_{k}, \tilde E_6, \tilde E_7, \tilde
E_8$-types and the corresponding integrable models should be
described explicitly. A systematic investigation of reductions of matrix
models constructed in this paper seems to be an interesting problem.
In this paper we have considered also integrable PDEs of
the chiral model type related to $PM$-structures with $m=2$. Compatible
associative multiplications described by $M$ and $PM$-structures can be
used for construction of integrable PDEs of the Landau-Lifshitz type in accordance with results
of \cite{golsok0}.

\vskip.3cm \noindent {\bf Acknowledgments.}
The second author (V.S.) is grateful to USTC (Hefei, China) for hospitality. The research was partially
supported by: RFBR grants 05-01-00775, 05-01-00189 and NSh grants 1716.2003.1,
2044.2003.2.


\vskip.3cm

\begin{thebibliography}{10}

\bibitem{golsok1} I.Z. Golubchik  and  V.V.  Sokolov,
\newblock{\it Compatible Lie brackets and integrable equations of the
principle chiral model type}, Func. Anal. and Appl., {\bf 36}(3),
172--181, 2002.

\bibitem{golsok2}  I.Z. Golubchik  and  V.V.  Sokolov,
 \newblock{\it Factorization of the loop algebras and compatible Lie
brackets}, Journal of Nonlinear Math. Phys., {\bf 12}(1), 343--350,
2005.

\bibitem{golsok3}  I.Z. Golubchik  and  V.V.  Sokolov,
 \newblock{\it Compatible Lie brackets and Yang-Baxter
equation}, Theoret. and Math. Phys.,  {\bf 146}(2),
 159--169, 2006.

 \bibitem{golsok4}  I.Z. Golubchik  and  V.V.  Sokolov,
 \newblock{\it Factorization of the loop algebra and integrable top-like
systems}, Theoret. and Math. Phys.,  {\bf 141}(1),
 3--23, 2004.

\bibitem{bolsbor} A.V. Bolsinov and A.V. Borisov,
\newblock{\it Lax representation and compatible Poisson brackets on Lie
algebras}, Math. Notes, {\bf 72}(1), 11-34, 2002.

\bibitem{odsok}  A.V. Odesskii  and V. V. Sokolov,   \newblock{\it Algebraic structures
connected with pairs of compatible associative algebras}, submitted
to IMRN,  (2006).


\bibitem{man} S.V.~Manakov, {\it Note on the integration of Euler's
equations of the dynamics of an n-dimensional rigid body},  Funct.
Anal. Appl. {\bf 10}, no.{\bf 4} (1976), 93-94

\bibitem{fed} Fedorov Yu. N. {\it Integrable systems, Lax
representations, and confocal quadrics. Dynamical systems in
classical mechanics},  Amer. Math. Soc., Providence, RI,
Amer. Math. Soc. Transl. Ser. 2, {\bf 168}, 173--199, 1995.

\bibitem{miksok} Mikhailov A.V., Sokolov V.V. {\it Integrable ODEs
on Associative Algebras},  Comm. in Math. Phys., {\bf 211}(1), 231--251, 2000.

\bibitem{pohl}  K. Pohlmeyer, {\it Integrable Hamiltonian systems and interactions through quadratic
constraints}, Commun. Math. Phys.  {\bf 46},  207--221, 1976.


\bibitem{stek1}  Stekloff W. {\it Sur le mouvement d\'un corps solide ayant
une cavite de forme ellipsoidale remple par un liquide
incompressible en sur les variations des latitudes},  Ann. Fac. Sci.
Toulouse Sci. Math. Sci. Phys. (3) 1 (1909), 145--256

\bibitem{puan} Poincare H. {\it Sur la precession des corps deformables},
Bull. Astr., 1910, v. 27, p. 321--356.

\bibitem{shot} Schottky F. {\it Uber das analytische Problem der Rotation
eines starren K\"orpers in Raume von vier Dimensionen},
Sitzungsberichte drer K\"onigligh preussischen Academie der
Wissenschaften zu Berlin, 1891, v. XIII, p. 227--232.

\bibitem{rat} Ratiu T. {\it Euler-Poisson equations on Lie
algebras and the $N$-dimensional heavy rigid body}, Amer. J. Math., {\bf 103}(3),
409--448. 1982.

\bibitem{bf}  Bolsinov A.V., Fedorov Yu. N. {\it Multi-dimensional
integrable generalization of the Steklov-Lyapunov systems}, Vestnik MGU, {\bf 6}, 53--56, 1992.


\bibitem{magri} F. Magri, \newblock{\it A simple model of the integrable
Hamiltonian equation}, J. Math. Phys., {\bf 19}, 1156--1162, 1978.


\bibitem{golsok0} I. Z. Golubchik  and  V. V.  Sokolov,
\newblock{\it Generalized Heizenberg equations on Z-graded Lie algebras},
Theoret. and Math. Phys, {\bf 120}(2), 1019--1025, 1999.





\end{thebibliography}
\end{document}